%% file: BMZ17.tex
\documentclass[graybox]{svmult}

\usepackage{mathptmx}       
\usepackage{helvet}         
\usepackage{courier}        
\usepackage{type1cm}        
\usepackage{makeidx}         
\usepackage{graphicx}        
\usepackage{multicol}        
\usepackage[bottom]{footmisc}

\usepackage{amsmath}
\usepackage{amsfonts}
\usepackage{xcolor}
\usepackage[utf8]{inputenc}

\usepackage{hyperref}
\hypersetup{colorlinks=true,linkcolor=blue,urlcolor=blue}

\def\mig{\frac{1}{2}}

\def\h2{\frac{h}{2}}
\def\u2{\frac{1}{2}}

\def\h{\tilde{h}}

\def\bigO{\mathcal{O}}

\makeindex             

\begin{document}

\title*{High order in space and time schemes through an approximate
  Lax-Wendroff procedure}
\titlerunning{An approximate Lax-Wendroff procedure}
\author{A.~Baeza, P.~Mulet and D.~Zorío}
\institute{
Departament de Matemàtiques,  Universitat de
  València    (Spain); emails: antonio.baeza@uv.es, mulet@uv.es, david.zorio@uv.es. This research was partially supported by Spanish MINECO grants
   MTM2011-22741 and MTM2014-54388.
}
%
%

\leavevmode\thispagestyle{empty}

\noindent This version of the article has been accepted for publication, after a peer-review process, and is subject to Springer Nature’s AM terms of use, but is not the Version of Record and does not reflect post-acceptance improvements, or any corrections. The Version of Record is available online at: \url{https://doi.org/10.1007/978-3-319-65870-4_31}

\newpage

\maketitle

\abstract{In this work we present an extension of the scheme developed
  by Qiu and Shu in 2003 to numerically solve hyperbolic conservation
  laws, with arbitrarily high 
  order both in space and time, based on the Lax-Wendroff procedure, using the 
  finite difference Shu-Osher's scheme and the WENO spatial
  reconstruction technique. 
  The aforementioned extension mainly consists in some improvements based on 
  simplifications in the computation of high order terms. On the one hand, a technique 
  which avoids the computation of the derivatives of the flux (which 
  commonly requires symbolic manipulation software, in addition to being 
  computationally expensive) is developed. On the other hand, a more careful discontinuity 
  detection is implemented in order to avoid a propagation of large terms at the 
  approximations of the derivatives.}

\section{Introduction}
\label{sec:intro}
In this paper we present a high order accurate temporal scheme, whose
derivation takes as starting point the one that was proposed in 2003
by Qiu and Shu \cite{QiuShu2003},
for numerically solving hyperbolic conservation laws,
based on the conversion of time derivatives 
to spatial derivatives of another expression involving
temporal derivatives of the flux, 
through the Cauchy-Kowalewski technique, following the Lax-Wendroff
procedure. 
The novelty in this work is the replacement of 
the exact derivatives of the flux by approximations of a suitable order in
order to reduce both the implementation and the computational cost 
maintaining the global order of the method, as
well as a fluctuation control method which avoids the expansion of large
terms at the discretization of the high order derivatives. The method is applied 
to a Shu-Osher finite-difference spatial discretization \cite{ShuOsher1989}.


The paper is organized as follows: in Section \ref{sec:ns} we
introduce the numerical scheme originally proposed by Qiu and Shu in
\cite{QiuShu2003}; in Section \ref{sec:cka}, we present the novel 
approach, based on using central difference approximations for the time
derivatives of the flux; Section \ref{sec:fc} stands for the
fluctuation control, a technique to avoid the excessive propagation of
diffusion around a discontinuity; in Section \ref{sec:ne} some
numerical experiments are presented; finally, some conclusions are
drawn in Section \ref{sec:conc}.

\section{Numerical scheme}\label{sec:ns}

Our goal is to numerically solve systems of $d$-dimensional $m$
hyperbolic conservation laws

\begin{equation*}\label{hcl}
    u_t+\sum_{i=1}^df^i(u)_{x_i}=0.
\end{equation*}

For the sake of simplicity, we start with the one-dimensional scalar
case ($d=m=1$). For the solution $u(x, t)$ of $u_t+f(u)_x=0$ on a fixed 
spatial   grid ($x_i$) with spacing $h=x_{i+1}-x_{i}$ and some time $t_n$
   from a temporal grid with spacing $\delta=\Delta t=t_{n+1}-t_n>0$,
   proportional to $h$, $\delta=\tau h$, where $\tau$ is dictated by
   stability restrictions (CFL condition) we use
the following  notation for time derivatives of $u$ and $f(u)$:
   \begin{align*}
     u_{i,n}^{(l)}&=\frac{\partial^{l} u(x_i, t_n)}{\partial t^l},\\
     f_{i,n}^{(l)}&=\frac{\partial^{l} f(u)(x_i, t_n)}{\partial t^l}.
   \end{align*}
Our goal is to obtain an $R$-th order accurate numerical scheme, i.e., a scheme
with a local truncation error 
of order $R+1$, based on the Taylor expansion  of the solution $u$
from time $t_n$ to time $t_{n+1}$:
$$u_i^{n+1}=\sum_{l=0}^R\frac{\Delta
  t^l}{l!}u_{i,n}^{(l)}+\bigO(\Delta t^{R+1}).$$
To achieve this we aim to define corresponding approximations 
\begin{align*}
  \widetilde u_{i,n}^{(l)}&=u_{i,n}^{(l)}+\bigO(h^{R+1-l}),\\
  \widetilde f_{i,n}^{(l)}&=f_{i,n}^{(l)}+\bigO(h^{R-l}),
\end{align*}
by recursion on $l$, assuming  (for a local truncation error
analysis) that  $\widetilde
u^{0}_{i,n}=u^{(0)}_{i,n}=u(x_i, t_n)$.

 The fact that $u$ solves the system of conservation laws implies that 
the time derivatives
$u_{i,n}^{(l)}$, $1\leq l\leq R$, can be written in terms of
the first spatial
derivative of some function of $u_{i,n}^{(j)}$, $j<l$, 
by using the chain rule on $f$, which can be written as 
\begin{equation}\label{eq:35}
f_{i,n}^{(l-1)}=F_{l-1}(u_i^n,u_{i,n}^{(1)},\ldots,u_{i,n}^{(l-1)})
\end{equation}
for some function $F_{l-1}$, and following the Cauchy-Kowalewski (or Lax-Wendroff for
second order) procedure:
\begin{equation}\label{eq:ck}
\frac{\partial^lu}{\partial t^l}=\frac{\partial^{l-1}}{\partial
  t^{l-1}}
\big(u_t\big)=-\frac{\partial^{l-1}}{\partial t^{l-1}}\big(f(u)_x\big)=
-\left[\frac{\partial^{l-1}f(u)}{\partial t^{l-1}}\right]_x,
\end{equation}

Specifically, to approximate the first time derivative, $u_t=-f(u)_x$,
we use the Shu-Osher finite difference scheme \cite{ShuOsher1989}
with upwinded WENO spatial
reconstructions  \cite{JiangShu1996} of order $2r-1$ in the flux function:
\begin{equation}\label{eq:40}
u_{i,n}^{(1)}=u_t(x_i,t_n)=-[f(u)]_x(x_i,t_n)
=-\frac{\hat{f}_{i+\frac{1}{2}}^n-\hat{f}_{i-\frac{1}{2}}^n}{h}+\mathcal{O}(h^{2r-1}).
\end{equation}
Much cheaper centered differences are used instead for the next derivatives.
We expound the general procedure for a third order accurate scheme
($R=3$) for a scalar  one-dimensional conservation law.
Assume we have numerical data,
$\{\widetilde u_i^n\}_{i=0}^{M-1}$, which approximates $u(\cdot, t_n)$ and
want to compute
an approximation for $u(\cdot, t_{n+1})$ at the same nodes, namely,
$\{\widetilde u_i^{n+1}\}_{i=0}^{M-1}$.

First, we compute an approximation of $u_t$ by the procedure stated
above: 
$$\widetilde
u_{i,n}^{(1)}=-\frac{\hat{f}_{i+\mig}^n-\hat{f}_{i-\mig}^n}{h},$$
with
$$\hat{f}_{i+\mig}^n=\hat{f}(\widetilde u_{i-r+1}^n,\ldots,\widetilde
u_{i+r}^n)$$
the numerical fluxes, which are obtained through upwind WENO spatial
reconstructions of order $2r-1$, with $r=\lceil\frac{R+1}{2}\rceil=2$.

Once the corresponding nodal data is obtained for the approximated
values of $u_t$, we compute
$$u_{tt}=[u_t]_t=[-f(u)_x]_t=-[f(u)_t]_x=-[f'(u)u_t]_x,$$
where $f'(u)u_t$ is now an approximately known expression for the
required nodes. We 
use then a second order centered difference in order to obtain the
approximation:
$$\widetilde u_{i,n}^{(2)}=-\frac{\widetilde
  f^{(1)}_{i+1,n}-\widetilde f^{(1)}_{i-1,n}}{2h},$$
where
$$\widetilde
f^{(1)}_{i,n}=F_1(\widetilde u_{i,n}^{(0)},\widetilde
u_{i,n}^{(1)})=f'(\widetilde u_{i,n}^{(0)})\widetilde u_{i,n}^{(1)},$$

Finally, we approximate the third derivative:
$$u_{ttt}=[u_t]_{tt}=[-f(u)_x]_{tt}=-[f(u)_{tt}]_x
=-\Big(f''(u)u_t^2+f'(u)u_{tt}\Big)_x,$$
where again the function $f''(u)u_t^2+f'(u)u_{tt}$ is approximately
known at the nodes and therefore $u_{ttt}$ can be
computed by second order accurate centered differences (note that 
in this case it would be required only a first order accurate
approximation; however, the order of centered approximations is always
even):
$$\widetilde u_{i,n}^{(3)}=-\frac{\widetilde
  f^{(2)}_{i+1,n}-\widetilde f^{(2)}_{i-1,n}}{2h},$$
where
$$\widetilde f^{(2)}_{i,n}=F_2(\widetilde u_{i,n}^{(0)},\widetilde
u_{i,n}^{(1)},\widetilde u_{i,n}^{(2)})
=f''(\widetilde u_{i,n}^{(0)})\cdot(\widetilde
u_{i,n}^{(1)})^2+f'(\widetilde u_{i,n}^{(0)})\cdot(\widetilde
u_{i,n}^{(2)})^2.$$
Once all the needed data is obtained, we advance in time by replacing
the terms of the
third order Taylor expansion in time of $u(\cdot, t_{n+1})$ by their
corresponding nodal approximations:
\begin{equation*}
  \begin{split}
    \widetilde u_i^{n+1}=\widetilde u_i^n+\Delta t\widetilde u_{i,n}^{(1)}+\frac{\Delta
      t^2}{2}\widetilde u_{i,n}^{(2)}+\frac{\Delta t^3}{6}\widetilde
    u_{i,n}^{(3)}.
  \end{split}
\end{equation*}

As we shall see, the above example can be extended to arbitrarily
high order time 
schemes through the computation of the suitable high order central
differences of the nodal values
$$\widetilde f^{(l)}_{i,n}=F_l(\widetilde u_{i,n}^{(0)},\widetilde
u_{i,n}^{(1)},\ldots,\widetilde u_{i,n}^{(l)})
=f_{i,n}^{(l)}+\mathcal{O}(h^{R-l+1}).$$
The generalization to multiple dimensions is straightforward, since 
now the Cauchy-Kowalewski procedure, being based on the fact that
$u_t=-\nabla\cdot f(u)$, yields
\begin{equation*}
\frac{\partial^l u}{\partial t^l}=-\nabla \cdot
\Big(\frac{\partial^{l-1}f(u) }{\partial t^{l-1}}\Big)
=-\sum_{i=1}^{d}\frac{\partial}{\partial x_i}\left(\frac{\partial^{l-1}f^{i}(u) }{\partial t^{l-1}}\right)
\end{equation*}
and that the spatial reconstruction procedures are done separately for
each dimension.  
For the case of the systems of equations, the time derivatives are now
computed through tensorial products of the corresponding derivatives
of the Jacobian of the fluxes.
The general procedure for systems and multiple
dimensions is thus easily generalizable and further details about the
procedure can be found in \cite{QiuShu2003, ZorioBaezaMulet}. 
Closed formulas to
explicitly compute the above expressions can be found in the
literature, such as the Fa\`a di Bruno formula \cite{FaaDiBruno1857}.

\section{The approximate Lax-Wendroff procedure}\label{sec:cka}

As reported by the authors of \cite{QiuShu2003}, the computation of
the exact nodal values of $f^{(k)}$ can be very
expensive as $k$ increases, since the number of  required operations
increases exponentially. Moreover, implementing it is costly and requires
large symbolic computations for each equation.
We now present an alternative, which is much less expensive for large $k$
and less dependent on the equation, in the sense that its only
requirement is the knowledge of the flux function.
This procedure also works in the multidimensional case and in the case
of systems as well (by
working componentwise). This
technique is based on the observation that  approximations
$\widetilde f^{(l-1)}\approx f^{{(l-1)}}$ can be easily obtained by
a suitable use of suitable finite differences, rather than using the
exact expression $F_{l-1}$ in \eqref{eq:35}.

Let us introduce some notation for the case of a
 one-dimensional system, that we assume for the sake of simplicity.
For a function $u\colon \mathbb R\to \mathbb R^{m}$, we denote the
discretization of the  
function on  the grid defined by a base
  point $a$ and grid space $h$ by
  \begin{equation*}
    G_{a,h}(u)\colon\mathbb{Z} \to \mathbb {R}^{m},\quad
    G_{a,h}(u)_i=u(a+ih).
  \end{equation*}
The symbol  $ \Delta^{p,q}_{h}$ denotes 
the centered finite differences operator that  approximates $p$-th order
  derivatives to order $2q$ on grids with spacing $h$. For any $u$
  sufficiently differentiable, it satisfies:
  \begin{align}\label{eq:3}
    \Delta^{p,q}_{h}
    G_{a,h}(u)=u^{(p)}(a)+\alpha^{p,q}u^{(p+2q)}(a)h^{2q}+\bigO(h^{2q+2}).
  \end{align}

We aim to define approximations $\widetilde u^{(k)}_{i,n} \approx
u^{(k)}_{i,n}$, $k=0,\dots,R$, recursively. We start the recursion with
\begin{equation}\label{eq:151}
  \begin{aligned}
    \widetilde u^{(0)}_{i,n}&=u_{i}^{n},\\
    \widetilde
    u^{(1)}_{i,n}&=-\frac{\hat{f}_{i+\frac{1}{2}}^n-\hat{f}_{i-\frac{1}{2}}^n}{h},
  \end{aligned}
\end{equation}   
  where $\hat f_{i+\mig}^{n}$ are computed from the known data $(u_{i}^{n})$
  by applying upwind WENO
  reconstructions (see
  \cite{ShuOsher1989,DonatMarquina96,JiangShu1996} for
  further details).

Associated to fixed $h, i, n$, once obtained $\widetilde
u^{(l)}_{i,n}$, $l=0,\dots,k$, in the recursive process we define
the $k$-th degree approximated Taylor polynomial $T_k[h,i, n]$ by
\begin{align*}
  T_k[h,i, n](\rho)=\sum_{l=0}^{k}\frac{\widetilde{u}^{(l)}_{i,n} }{l!} \rho^l.
\end{align*}

By recursion,  for $k=1,\dots,R-1$, we define
\begin{equation}\label{eq:15}
  \begin{aligned}
    \widetilde f^{(k)}_{i,n} &=
    \Delta_{\delta}^{k,\left\lceil \frac{R-k}{2}\right\rceil}\Big(G_{0,\delta}\big(f(T_k[h, i,
    n])\big)\Big),\\
    \widetilde
    u^{(k+1)}_{i,n}&=-\Delta_h^{1,\left\lceil\frac{R-k}{2}\right\rceil}
    \widetilde f^{(k)}_{i+\cdot,n},
  \end{aligned}
\end{equation}   
where we denote by $\widetilde
f^{(k)}_{i+\cdot,n}$ the vector 
given by the elements $(\widetilde f^{(k)}_{i+\cdot,n})_{j}=\widetilde
f^{(k)}_{i+j,n}$ and $\delta=\Delta t $. With all these ingredients, the
proposed scheme is:
\begin{equation}\label{eq:60}
u_i^{n+1}=u_{i}^{n}+\sum_{l=1}^R\frac{\Delta
  t^l}{l!}\widetilde u_{i,n}^{(l)}.
\end{equation}

It can be proven that the method resulting from this construction 
is $R-th$ order accurate and can be written in conservation form, 
see \cite{ZorioBaezaMulet}

\section{Fluctuation control}\label{sec:fc}
Now we focus on the computation of the approximate nodal values of the
first order time derivative. Typically, one would simply take the
approximations obtained through the upwinded reconstruction procedure
in the Shu-Osher's finite difference approach, that is,
\begin{equation}\label{eq:61}
\widetilde u_{j,n}^{(1)}=-\frac{\hat{f}_{j+\frac{1}{2},n}-\hat{f}_{j-\frac{1}{2},n}}{h}.
\end{equation}

However,taking directly these values as the first derivative used to compute
the next derivatives through \eqref{eq:60} can produce wrong results 
if the data is not smooth, as routinely happens in hyperbolic systems. 
In fact, it will include $\mathcal{O}(h^{-1})$ terms wherever there
is a discontinuity, which we will call from now on
\textit{fluctuations}. These terms will appear provided
$\hat{f}_{j-\frac{1}{2},n}$ and $\hat{f}_{j+\frac{1}{2},n}$ come from
different sides of a discontinuity (or some of them has mixed
information of both sides due to a previous flux splitting procedure to
reconstruct the interface values), since in that case
$\hat{f}_{j+\frac{1}{2},n}-\hat{f}_{j-\frac{1}{2},n}=\mathcal{O}(1)$.

In practice, this implies that the $k$-th derivative, $1\leq k\leq R$,
will have terms of magnitude $\mathcal{O}(h^{-k})$, therefore, the
term which appears on the Taylor expansion term, which is multiplied
by $\frac{\Delta t^k}{k!}$, a term of magnitude $\mathcal{O}(h^k)$,
will be ultimatelly $\mathcal{O}(1)$. This may result in undesired
diffusion, oscillations or even a complete failure of the scheme in
some cases.

Our proposal is to compute an alternative approximation of 
$\widetilde u_{j,n}^{(1)}$ as described in Section \ref{sec:cweno} and
replace \eqref{eq:61} by this new approximation for the recursive 
computation in \eqref{eq:15} only, 
maintaining \eqref{eq:61} for its use in \eqref{eq:60}, thus ensuring 
the proper upwindind.

\subsection{Central WENO reconstructions}\label{sec:cweno}

Let us assume that our spatial scheme is $(2r-1)$-th order accurate
WENO. After all the operations performed for the reconstruction of the
numerical fluxes at the interfaces, the stencil of points that is used in order to approximate
the derivative at the node $x_i$ is the following set of $2r+1$
points:
\begin{equation}
\label{eq:250}\{x_{i-r},\ldots,x_i,\ldots,x_{i+r}\},
\end{equation}
whose corresponding flux values, $f_j=f(u_j)$, are
$$\{f_{i-r},\ldots,f_i,\ldots,f_{i+r}\}.$$
The procedure that we next expound only uses information from the
stencil
\begin{equation}\label{eq:251}
  \mathcal{S}_{i+r-1}^{2r-1}:=\{i-r+1,\ldots,i,\ldots,i+r-1\},
\end{equation}  
thus ignoring  the flux values $f_{i-r}, f_{i+r}$ at the edges of the
stencil in \eqref{eq:250}.

For fixed $i$, let $q_k^{r}$ be the interpolating polynomial of
degree $\leq r-1$ such that $q_k^{r}(x_j)=f_j,$
$j\in\mathcal{S}_{i+k}^{r}:=\{i+k-r+1,\dots, i+k,\}$, $0\leq k\leq r-1$.
After the previous discussion, our goal is to obtain an approximation
of the flux derivative $f(u)_{x}(x_i)$ from the stencil
$\mathcal{S}_{i+r-1}^{2r-1}$ which is $(2r-1)$-th order accurate  if the
nodes in the stencil lie within a smoothness region for $u$ or is
$\bigO(1)$ otherwise. We use Weighted Essentially Non Oscillatory
techniques to achieve this purpose.

The following lemma is easily established.

\begin{lemma}
  There exists a set of constants $\{c^{r}_k\}_{k=1}^r$
satisfying $0<c^{r}_k<1$, for $0\leq k\leq r-1$,
$\sum_{k=0}^{r-1}c^{r}_k=1$,  such that 
$$\sum_{k=0}^{r-1}c^{r}_k(q^{r}_k)'(x_i)=(q_{r-1}^{2r-1})'(x_i).$$
\end{lemma}

If $f_j=f(u(x_j, t_n))$, for smooth enough $u$ and fix $t_n$, then
\begin{align}
  \label{eq:201}
  (q_{k}^{r})'(x_i)&=f(u)_x(x_i,  t_n)+d_{k}^{r}(x_i)h^{r-1}+\bigO(h^{r}), k=0,\dots,r-1,\\
  \label{eq:202}
  (q_{r-1}^{2r-1})'(x_i)&=f(u)_x(x_i,  t_n)+d_{r-1}^{2r-1}(x_i)h^{2r-2}+\bigO(h^{2r-1}).
\end{align}
for continuously differentiable $d_{k}^{r}, d_{r-1}^{2r-1}$. The  goal
is to obtain the accuracy in \eqref{eq:202} by  a suitable nonlinear convex
combination of \eqref{eq:201}
\begin{equation}\label{eq:203}
  \sum_{k=0}^{r-1}w^{r}_k(q^{r}_k)'(x_i)=f(u)_x(x_i,
  t)+\widetilde d_{r-1}^{2r-1}(x_i)h^{2r-2}+\bigO(h^{2r-1}),
\end{equation}  
where $w^{r}_k=c^{r}_{k}(1+\bigO(h^{r-1}))$ if the whole stencil
  $x_{i-r+1},\dots,x_{i+r-1}$ lies within a smoothness region for $u$
  and $w^{r}_{k}=\bigO(h^{r-1})$ if the $k$-th stencil crosses a
  discontinuity and there are at least another stencil which does
  not. We follow  Weighted Essentially Non Oscillatory classical
  techniques  \cite{LiuOsherChan1994,JiangShu1996}. From now
  on we drop the superscript $r$ in $q_k^{r}$.
  
  Furthermore, we need the approximation in \eqref{eq:203} to be in
  conservation form. To achieve this we use 
  the polynomial $p_k$ of degree $r-1$ satisfying
$$\frac{1}{h}\int_{x_{j-\frac{1}{2}}}^{x_{j+\frac{1}{2}}}p_k(x)dx=f_j,\quad
i-r+1+k\leq j\leq i+k,\quad 0\leq k\leq r-1,$$ and $\widetilde p_k(x)$ a
primitive of it. It can be seen that the polynomial
\begin{equation*}
\widetilde {q}_k(x) = \frac{\widetilde p_{k}(x+\frac{h}{2})-\widetilde
  p_{k}(x-\frac{h}{2})}{h},
\end{equation*}
has degree $\leq r-1$ and that $\widetilde{q}_{k}(x_{j})=f_j$,
$j=i-r+1+k,\dots, i+k$, and therefore $\widetilde {q}_k(x)$ must
coincide with $q_{k}$. Thus
$$q_k'(x_j)=\frac{(\widetilde p_k)'(x_{j+\frac{1}{2}})-(\widetilde p_k)'(x_{j-\frac{1}{2}})}{h}=\frac{p_k(x_{j+\frac{1}{2}})-p_k(x_{j-\frac{1}{2}})}{h}.$$

Now, let us define the following Jiang-Shu smoothness indicators using
the definition of $p_k$:
\begin{equation}\label{smthi}
I_k=\sum_{\ell=1}^{r-1}\int_{x_{i-\frac{1}{2}}}^{x_{i+\frac{1}{2}}}
h^{2\ell-1}p_k^{(\ell)}(x)^2dx,\quad 0\leq k\leq r-1
\end{equation}
so that we can define the weights as follows:
\begin{equation}\label{eq:204}
  \begin{aligned}
\omega_k&=\frac{\alpha_k}{\sum_{l=1}^r\alpha_l},\quad \alpha_k=\frac{c_k}{(I_k+\varepsilon)^m},
\end{aligned}
\end{equation}
with $\varepsilon>0$ a small positive quantity, possibly depending on
$h$.
Following the techniques in \cite{SINUM2011}, since
$p_{k}^{(\ell)}-p_{j}^{(\ell)}=\bigO(h^{r-\ell})$ at regions of
smoothness, whereas $\int_{x_{i-\mig}}^{x_{i+\mig}}
(p_{k}')^2dx=\bigO(h^{-2})$ if the corresponding stencil
$\mathcal{S}_{k}^{r}$ crosses a discontinuity, the smoothness indicators satisfy
$I_{k}-I_{j}=\bigO(h^{r+1})$ at regions of smoothness and $I_{k}\not \to
0$ if the $k$-th stencil crosses a discontinuity. Therefore, the definition
\eqref{eq:204} satisfies the requirements mentioned above in order to
achieve maximal order even at  smooth extrema, provided that 
the parameter  $\varepsilon>0$, besides avoiding divisions by zero, is
chosen as $\varepsilon=\lambda h^2,$ with $\lambda\sim 
f(u)_x$, and that the exponent $m$ in \eqref{eq:204}  makes the weight
$\omega_k=\mathcal{O}(h^{r-1})$ wherever
there is a discontinuity at that stencil. Since one wants to attain
the maximal  possible order in such case, which corresponds to the
value interpolated from a smooth substencil, which is
$\mathcal{O}(h^r)$, then it suffices to set
$m=\lceil\frac{r}{2}\rceil$. Finally, we define $\widetilde{\widetilde{u}}_{i,n}^{(1)}$, the
smoothened approximation
of $u_t(x_i, t_n)$ that replaces $\widetilde{u}_{i,n}^{(1)}$ in \eqref{eq:15}
 as the result of the following convex combination:
$$\widetilde{\widetilde{u}}_{i,n}^{(1)}=-\sum_{k=1}^r\omega_kq'_k(x_i).$$

\section{Numerical experiments}\label{sec:ne}
In this section we present some 2D experiments with Euler equations
involving comparisons of the fifth order both in
space ($r=3$) and time ($R=2r-1=5$) exact and approximate Lax-Wendroff
schemes, together with the results obtained using the third order TVD
Runge-Kutta time discretization. 

From now on we will refer as WENO[]-LW[] to the exact Lax-Wendroff
procedure, WENO[]-LWA[] to the approximate Lax-Wendroff procedure,
WENO[]-LWF[] if a fluctuation control is used in the exact procedure,
WENO[]-LWAF[] if the fluctuation control comes together with the
approximate procedure and
WENO[]-RK[] when a Runge-Kutta method is used. In each case,
the first bracket stands for the value of the spatial accuracy order and
the second one for the time accuracy order.

\subsection{Smooth solution}

In order to test the accuracy of our scheme
in the general scenario of a multidimensional system of conservation
laws, we perform a test using the 2D Euler equations with smooth
initial conditions, given by
\begin{equation*}
  \begin{split}
    u_0(x,y)&=(\rho(x,y),v^x(x,y),v^y(x,y),E(x,y))\\
    &=\left(\frac{3}{4}+\frac{1}{2}\cos(\pi(x+y)),
    \frac{1}{4}+\frac{1}{2}\cos(\pi(x+y)),\right.\\
    &\left.\frac{1}{4}+\frac{1}{2}\sin(\pi(x+y)),
    \frac{3}{4}+\frac{1}{2}\sin(\pi(x+y))\right),
  \end{split}
\end{equation*}
where $x\in\Omega=[-1,1]\times[-1,1]$, with periodic boundary
conditions.

In order to perform the smoothness analysis, we compute a reference
solution in a fine mesh and then compute numerical solutions for the
resolutions $n\times n$, for $n=10\cdot 2^k,$
$1\leq k\leq 5$, obtaining the results shown in Tables
\ref{SE2-CK}-\ref{SE2-CKRS} at the time $t=0.025$ for $\textnormal{CFL}=0.5$.
\begin{table}[htb]
  \centering
  \begin{tabular}{|c|c|c|c|c|}
    \hline
    $n$ & Error $\|\cdot\|_1$ & Order $\|\cdot\|_1$ & Error
    $\|\cdot\|_{\infty}$ & Order $\|\cdot\|_{\infty}$ \\
    \hline
    40 & 1.80E$-5$ & $-$ & 2.74E$-4$ & $-$ \\
    \hline
    80 & 1.09E$-6$ & 4.05 & 1.80E$-5$ & 3.93 \\
    \hline
    160 & 3.89E$-8$ & 4.80 & 7.36E$-7$ & 4.61 \\
    \hline
    320 & 1.29E$-9$ & 4.92 & 2.49E$-8$ & 4.88 \\
    \hline
    640 & 4.11E$-11$ & 4.97 & 8.07E$-10$ & 4.95 \\
    \hline
    1280 & 1.23E$-12$ & 5.06 & 2.43E$-11$ & 5.06  \\
    \hline
  \end{tabular}
  \caption{Error table for 2D Euler equation, $t=0.025$. WENO5-LW5.}
  \label{SE2-CK}
\end{table}
\begin{table}[htb]
  \centering
  \begin{tabular}{|c|c|c|c|c|}
    \hline
    $n$ & Error $\|\cdot\|_1$ & Order $\|\cdot\|_1$ & Error
    $\|\cdot\|_{\infty}$ & Order $\|\cdot\|_{\infty}$ \\
    \hline
    40 & 1.80E$-5$ & $-$ & 2.74E$-4$ & $-$ \\
    \hline
    80 & 1.09E$-6$ & 4.05 & 1.80E$-5$ & 3.93 \\
    \hline
    160 & 3.89E$-8$ & 4.80 & 7.36E$-7$ & 4.61 \\
    \hline
    320 & 1.29E$-9$ & 4.92 & 2.49E$-8$ & 4.88 \\
    \hline
    640 & 4.11E$-11$ & 4.97 & 8.07E$-10$ & 4.95 \\
    \hline
    1280 & 1.23E$-12$ & 5.06 & 2.43E$-11$ & 5.06  \\
    \hline
  \end{tabular}
  \caption{Error table for 2D Euler equation, $t=0.025$. WENO5-LWA5.}
  \label{SE2-CKR}
\end{table}
\begin{table}[htb]
  \centering
  \begin{tabular}{|c|c|c|c|c|}
    \hline
    $n$ & Error $\|\cdot\|_1$ & Order $\|\cdot\|_1$ & Error
    $\|\cdot\|_{\infty}$ & Order $\|\cdot\|_{\infty}$ \\
    \hline
    40 & 2.63E$-5$ & $-$ & 2.97E$-4$ & $-$ \\
    \hline
    80 & 1.58E$-6$ & 4.06 & 2.01E$-5$ & 3.89 \\
    \hline
    160 & 6.66E$-8$ & 4.57 & 1.06E$-6$ & 4.24 \\
    \hline
    320 & 2.33E$-9$ & 4.84 & 4.08E$-8$ & 4.70 \\
    \hline
    640 & 7.60E$-11$ & 4.94 & 1.34E$-9$ & 4.93 \\
    \hline
    1280 & 2.35E$-12$ & 5.02 & 4.06E$-11$ & 5.04  \\
    \hline
  \end{tabular}
  \caption{Error table for 2D Euler equation, $t=0.025$. WENO5-LWAF5.}
  \label{SE2-CKRS}
\end{table}
We can see thus that our scheme achieves the desired accuracy in
the general scenario of a multidimensional system of conservation
laws. Also, we can
see that the results obtained through the approximate Lax-Wendroff
procedure are almost the same than those obtained using the exact
version. The version with fluctuation control also yields the desired
accuracy.

\subsection{Double Mach Reflection}

This experiment uses the Euler equations to model a vertical right-going Mach
10 shock colliding with an equilateral triangle. By symmetry, this is
equivalent to a collision with a ramp with a slope of 30 degrees with
respect to the horizontal line.

For the sake of simplicity, we consider the equivalent problem of 
an oblique shock whose vertical angle is
$\frac{\pi}{6}$ rad in the rectangle $\Omega=[0,4]\times[0,1].$  
 The initial conditions of the problem 
are
$$u_0(x,y)=\begin{cases}
  C_1 & y\leq\frac{1}{4}+\tan(\frac{\pi}{6})x,\\
  C_2 & y>\frac{1}{4}+\tan(\frac{\pi}{6})x,\\
\end{cases}$$
where
\begin{equation*}
  \begin{split}
    C_1=(\rho_1,v^x_1,v^y_1,E_1)^T
    &=(8,8.25\cos(\frac{\pi}{6}),-8.25\sin(\frac{\pi}{6}),563.5)^T,\\
    C_2=(\rho_2,v^x_2,v^y_2,E_2)^T
    &=(1.4,0,0,2.5)^T.\\
  \end{split}
\end{equation*}
We impose inflow boundary conditions, with value $C_1$, at the left
side, $\{0\}\times[0,1]$, outflow boundary conditions both at
$[0,\frac{1}{4}]\times\{0\}$ and $\{4\}\times[0,1]$, reflecting
boundary conditions at  $]\frac{1}{4},4]\times\{0\}$ and inflow
boundary conditions at the upper side, $[0,4]\times\{1\}$, which
mimics the shock at its actual traveling speed:
$$u(x,1,t)=\begin{cases}
  C_1 & x\leq\frac{1}{4}+\frac{1+20t}{\sqrt{3}}, \\
  C_2 & x>\frac{1}{4}+\frac{1+20t}{\sqrt{3}}.\\
\end{cases}$$
We run different simulations until $t=0.2$ at a resolution of
$2048\times512$ points
for $\textnormal{CFL}=0.4$ and a different combination of techniques, involving
WENO5-RK3, WENO5-LW5 and WENO5-LWA5.

\begin{figure}[htb]
    \centering
  \begin{tabular}{cc}
    \includegraphics[width=0.47\textwidth]{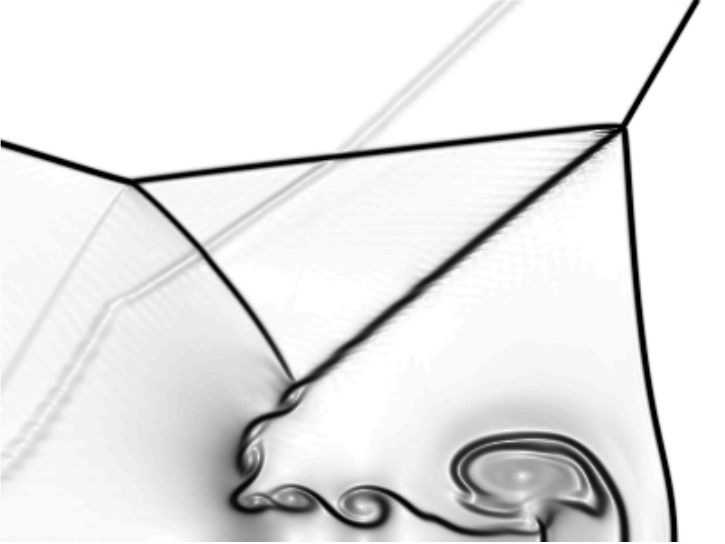}
    & \includegraphics[width=0.47\textwidth]{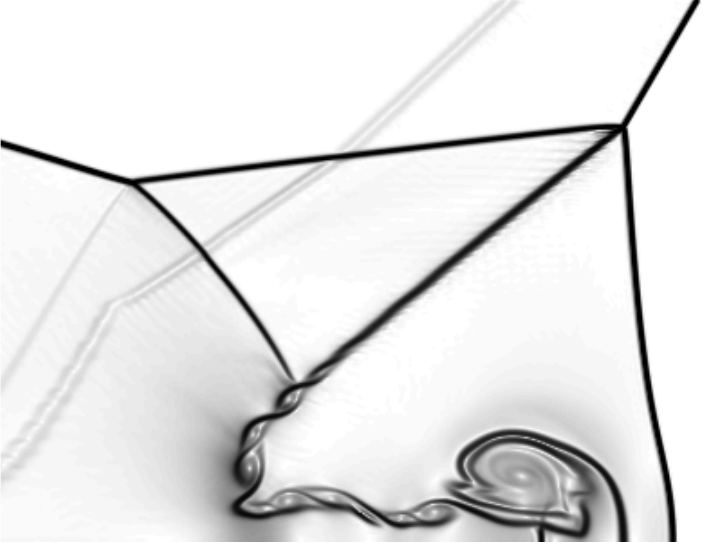} \\
    (a) WENO5-RK3 & (b) WENO5-LW5 \\[1ex]
    \includegraphics[width=0.47\textwidth]{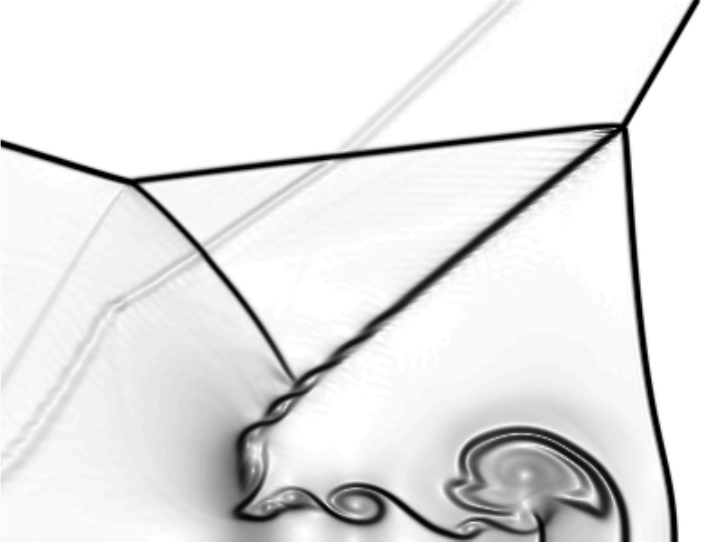}
    & \includegraphics[width=0.47\textwidth]{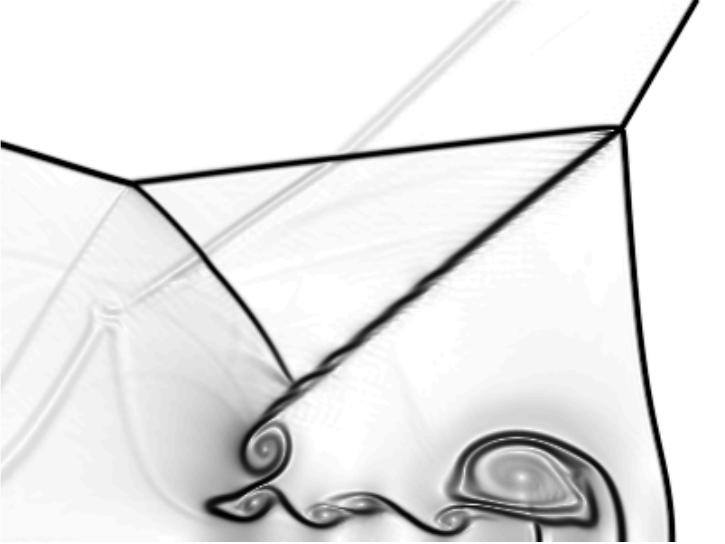} \\
    (c) WENO5-LWA5 & (d) WENO5-LWAF5 \\
  \end{tabular}
  \caption{Double Mach Reflection results. Density field.}
  \label{dmr}
\end{figure}

The results are presented in Figure \ref{dmr} as Schlieren plots of the turbulence zone.
It can be concluded that the results obtained
through the exact and approximate Lax-Wendroff techniques are
quite similar, and that the results obtained through the technique
with fluctuation control provides a slightly sharper profile.

Finally, in order to illustrate that the LW techniques are more
efficient than the RK time discretization, we show a performance test
involving the computational time required by each technique by running
the Double Mach Reflection problem for the resolution
$200\times50$. The results are shown in Table \ref{performance}, where
the field ``Efficiency'' stands for
$\displaystyle\frac{t_{\text{RK3}}}{t_{\text{LW*}}}$.
\begin{table}
  \centering
  \begin{tabular}{|c|c|}
    \hline
    Method & Efficiency \\
    \hline
    WENO5-LW5 & $1.44$ \\
    \hline
    WENO5-LWA5 & $1.54$ \\
    \hline
    WENO5-LWF5 & $1.33$ \\
    \hline
    WENO5-LWAF5 & $1.44$ \\
    \hline
  \end{tabular}
  \caption{Performance table.}
  \label{performance}
\end{table}

We can see from Table \ref{performance} that even the fifth order
Lax-Wendroff technique is more efficient than the third order
accurate Runge-Kutta scheme.

On the other hand, we see that the version with approximate fluxes has
a better performance than the main formulation, since less
computations are required for high order derivatives. On the other
hand, if the fluctuation control is used then the performance is
lower, however, the combination of the approximate
fluxes with the fluctuation control yields a fifth order accurate with
approximately the same efficiency than the original formulation, but
providing better results.

\section{Conclusions}\label{sec:conc}

In this paper we have developed an arbitrarly high order
Lax-Wendroff-type time schemes which does not require symbolic
computations to implement them, unlike those proposed by Qiu and Shu
in \cite{QiuShu2003}. Moreover, a fluctuation control has been
developed in order to avoid the propagation of incorrect data around
the discontinuities.

The results obtained in the numerical experiments are satisfactory,
and show that the approximate procedure yields essentially the same
results than the exact version with a much lower implementation cost
and less computationally expensive. On the other hand, the version with
fluctuation control, albeit increasing the computational cost,
produces numerical solutions with better resolution.

\input{referenc}

\end{document}

%% file: referenc.tex
%
%
%